\documentclass[12pt,a4paper]{article}
\usepackage{amsmath}
\usepackage{amsthm}
\usepackage{amsfonts}
\usepackage{amssymb}
\usepackage{euscript}
\usepackage[cp1251]{inputenc}
\usepackage[english] {babel}
\usepackage[pdftex]{graphicx}
\usepackage[10pt]{extsizes}

\begin{document}
\let\thefootnote\relax\footnote{This work was supported by the program “Leading Scientific Schools” under grant NSh-4833.2014.1.}

\begin{center}
  \textbf{Explicit characterization of some commuting differential operators of rank 2}
\end{center}
\begin{center}
  \textbf{V.Oganesyan}
\end{center}
\begin{center}
\textbf{1.Introduction}
\end{center}
If two differential operators
\begin{equation*}
L_n= \sum\limits^{n}_{i=0} u_i(x)\partial_x^i,  \quad  L_m= \sum\limits^{m}_{i=0} v_i(x)\partial_x^i
\end{equation*}
commute, then there is a nonzero polynomial $R(z,w)$ such that  $R(L_n,L_m)=0$ (see ~\cite{Chaundy}). The curve $\Gamma$ defined by $R(z,w)=0$ is called the \emph{spectral curve}. The genus of the curve $R(z,w)=0$ is called the genus of commuting pair. If\\
\begin{equation*}
L_n \psi=z\psi, \quad  L_m \psi=w\psi,
\end{equation*}
then $(z,w) \in \Gamma$. For almost all $(z,w) \in \Gamma$ the dimension of the space of common eigenfunctions $\psi$ is the same. The dimension of the space of common eigenfunctions of two commuting differential operators is called the \emph{rank}. The rank is a common divisor of m and n.\\
If the rank equals 1, then there are explicit formulas for coefficients of commutative operators in terms of Riemann theta-functions (see ~\cite{theta}).\\
The case when rank is greater than one is much more difficult. The first examples of commuting ordinary scalar differential operators of the nontrivial ranks 2 and 3 and the nontrivial genus g=1 were constructed by Dixmier ~\cite{Dixmier} for the nonsingular elliptic spectral curve $w^2=z^3-\alpha$, where $\alpha$ is an arbitrary nonzero constant:
\begin{equation*}
\begin{gathered}
L= (\partial_x^2 + x^3 + \alpha)^2 + 2x ,\\
M= (\partial_x^2 + x^3 + \alpha)^3 + 3x\partial_x^2 + 3\partial_x + 3x(x^2+\alpha) ,
\end{gathered}
\end{equation*}
where $L$ and $M$ is the commuting pair of the Dixmier operators of rank 2, genus 1. There is an example
\begin{equation*}
\begin{gathered}
L= (\partial_x^3 + x^2 + \alpha)^2 + 2\partial_x  ,\\
M= (\partial_x^3 + x^2 + \alpha)^3 + 3\partial_x^4 +  3(x^2+\alpha)\partial_x + 3x  ,
\end{gathered}
\end{equation*}
where $L$ and $M$ is the commuting pair of the Dixmier operators of rank 3, genus 1.\\
The general classification of commuting ordinary differential operators of rank greater than 1 was obtained by Krichever ~\cite{ringkrichever}. The general form of commuting operators of rank 2 for an arbitrary elliptic spectral curve was found by Krichever and Novikov ~\cite{novkrich}. The general form of operators of rank 3 for an arbitrary elliptic spectral curve (the general commuting operators of rank 3, genus 1 are parameterized by two arbitrary functions) was found by Mokhov ~\cite{Mokhov1},~\cite{Mokhov2}. Also operators of rank greater than 1 was considered in ~\cite{Grunbaum}, ~\cite{PrevWils}.\\
Mironov in ~\cite{Mironov} constructed examples of operators
\begin{equation*}
\begin{gathered}
L = (\partial_x^2 + A_3x^3+ A_2x^2 + A_1x + A_0 )^2 + g(g+1)A_3x ,\\
M^2 = L^{2g+1} + a_{2g}L^{2g} + ... + a_1L + a_0 ,
\end{gathered}
\end{equation*}
where $a_i$  are some constants and $A_i$, $A_3 \neq 0$, are arbitrary constants. Operators $L$ and $M$ are commuting operators of rank 2, genus g.
Furthermore, in ~\cite {Mironov2} it was proved by Mironov that  operators $L_1$ and $M_1$
\begin{equation*}
\begin{gathered}
L_1 = (\partial_x^2 + \alpha_1 \mathcal{P}(x) + \alpha_0)^2 + \alpha_1g_2g(g+1)\mathcal{P}(x), \quad \alpha_1 \neq 0 ,\\
M_1^2 = L_1^{2g+1} + a_{2g}L_1^{2g} + ... + a_1L_1 + a_0 ,
\end{gathered}
\end{equation*}
where $a_i$  are some constants, $\alpha_i$ are arbitrary constants and $\mathcal{P}$ satisfies the equation
\begin{equation*}
(\mathcal{P'}(x))^2 = g_2\mathcal{P}^2(x) + g_1\mathcal{P}(x) + g_0,  \quad g_2 \neq 0,
\end{equation*}
where $g_1$ and $g_2$ are arbitrary constants, is a commuting pair of rank 2, genus $g$.\\
Let $\wp(x)$ be the Weierstrass elliptic function satisfying the equation $(\wp'(x))^2 = 4\wp^3(x) + g_2\wp(x) + g_3$. Mironov proved in ~\cite {Mironov2} that  operators $L_2$ and $M_2$
\begin{equation*}
\begin{gathered}
L_2 = (\partial_x^2 + \alpha_1 \wp(x) + \alpha_0)^2 + s_1 \wp(x) + s_2 \wp^2(x),\\
M_2^2 = L_2^{2g+1} + b_{2g}L_2^{2g} + ... + b_1L_2 + b_0,
\end{gathered}
\end{equation*}
where $b_i$ are some constants, $\alpha_1=\frac{1}{4} - 2g^2 - 2g$, $s_1=\frac{1}{4}g(g+1)(16\alpha_0 + 5g_2)$,  $s_2=-4g(g+2)(g^2-1)$ and $\alpha_0$ is an arbitrary constant, are also a commuting pair of rank 2, genus $g$. Using the same methods many other examples were found  \cite{Davl1}, \cite{Davl2}.\\
Examples of commuting ordinary differential operators of arbitrary genus and arbitrary rank with polynomial coefficients were constructed in  ~\cite{Mokhov4} by Mokhov.\\
Consider
\begin{equation*}
\begin{gathered}
L_1=(\partial_x^2 + A_6x^6 + A_2x^2 )^2 + 16g(g+1)A_6x^4,\\
L_2=(\partial_x^2 + A_4x^4 + A_2x^2 + A_0)^2 + 4g(g+1)A_4x^2,
\end{gathered}
\end{equation*}
where $g \in \mathbb{N}$, $A_6 \neq 0$, $A_4\neq 0$, $A_2, A_0$ are arbitrary constants. Operators $L_1$ and $L_2$ commute with operators  of order 4g+2 (see~\cite{Vartan}, \cite{Vartan2}). The spectral curves of these operators have the form $w^2=z^{2m+1} + a_{2m}z^{2m} + ...+a_1z + a_0$. Moreover, the following theorems are proved in ~\cite{Vartan}, \cite{Vartan2}.
\\
\begin{itemize}
\item[1)]  If $L=(\partial_x^2 + A_nx^n + A_{n-1}x^{n-1} + ... + A_0)^2 + B_kx^{k} + B_{k-1}x^{k-1} + ...+B_0$, where $n>3, n \in \mathbb{N}$, $A_n \neq 0$, $B_k \neq 0$, commutes with a differential operator M of order $4g+2$ and M, L are operators of rank 2, then $k=n-2$ and $B_k=(n-2)^2m(m+1)A_n$ for some $m \in \mathbb{N}$.\\
\item[2)] For $n>6$  operator $L=(\partial_x^2 + A_nx^n)^2 + B_{n-2}x^{n-2}$ does not commute with any differential operator $M$ of order $4g+2$, where $M$ and $L$ could be a pair of rank 2.\\
\item[3)]  If $n=5$, then $L=(\partial_x^2 + Ax^5)^2 + 18Ax^{3}$, $A\neq 0$,  commutes with a differential operator M of order $4g+2$ for all $g$ and M,L are operators of rank 2. \\
\item[4)] The operator  $(\partial_x^2 + Ax^5)^2 + 9m(m+1)Ax^{3}$, $A\neq 0$ does not commute with any differential operator M of order $4g+2,$ for $m>1$, where M and L could be a pair of rank 2 .
\end{itemize}
Let us consider the operator
\begin{equation*}
L_4 = \partial_x^4 + u(x).
\end{equation*}
Assume that $u(x)$  has poles at points $a_1,a_2,...$. In a neighborhood of $a_i$
\begin{equation*}
u(x)=\frac{\varphi_{i,-k}}{(x-a_i)^k} + \frac{\varphi_{i,-k+1}}{(x-a_i)^{k-1}} +... +\varphi_{i,0} + \varphi_{i,1}(x-a_i)+ O((x-a_i)^2).
\end{equation*}
The main results of this paper are the following.\\
\\
\textbf{Theorem 1.1}\\
\emph{If $L_4=\partial_x^4 + u(x)$ commutes with a differential operator $M$ of order $4g+2$ and $M$, $L_4$ are operators of rank 2, then $u(x)$ can have pole only of order $4$, $\varphi_{i,-4}=n_i(4n_i+1)(4n_i+3)(4n_i+4)$, $n_i \in \mathbb{N}$, $\varphi_{i,4k-l}=0$, where $k=0,...,n_i$, $l=1,2,3$. Moreover  $\varphi_{i, 4r-1}=\varphi_{i, 4r-3}=0$, where $r=n_i+1,...,g$. Function $u(x)$ can't have isolated pole at infinity.}\\
\\
\textbf{Corollary 1.2}\\
\emph{Suppose $u(x)$ is elliptic, simply-periodic  or rational function and hasn't isolated singularity at infinity. Let $a_1$ be the unique pole in a fundamental parallelogram, a period-strip or the complex plane respectively. The operator $L_4=\partial^4_x + u(x)$ commutes with an operator $M$ of order $4g+2$ and $M$, $L_4$ are operators of rank 2 if and only if $\varphi_{1,-4}=g(4g+1)(4g+3)(4g+4)$, $\varphi_{1, 4k-l}=0$, where $k=0,...,g$, $l=1,2,3$.}\\
\\
\textbf{Corollary 1.3}\\
\emph{Let $\wp(x)$ be the Weierstrass elliptic function satisfying the equation $(\wp'(x))^2 = 4\wp^3(x) + g_2\wp(x) + g_3$. The Operator
\begin{equation*}
L_4=\partial^4_x + n(4n+1)(4n+3)(4n+4)\wp^2(x),
\end{equation*}
where  $n \in \mathbb{N}$, commutes with an operator of order $4n+2$ if and only if $\wp(x)$ is a solution of the equation $(\wp'(x))^2=4(\wp(x))^3 + g_2\wp(x)$. If $g_2=0$, then we obtain $L_4=\partial^4_x + \frac{n(4n+1)(4n+3)(4n+4)}{x^4}$. Calculations show that for $n$ less than $8$ the spectral curve is non-singular for almost all $g_2$.\\}
\\
\textbf{Example 1.4}\\
\emph{Let $\wp(x)$ be the Weierstrass elliptic function satisfying the equation $(\wp'(x))^2 = 4\wp^3(x) + g_2\wp(x) + g_3$. Consider the operator \begin{equation*}
L_4=\partial^4_x + 280\wp^2(x) + 280\wp^2(x-a),
\end{equation*}
Suppose that $g_3 = 0$ and $a=\omega_1$, $\omega_2$ or $\omega_1 + \omega_2$, where $\omega_i$ are half-periods. Then $L_4$ commutes with an operator of order $6$.}\\
\\
Analyzing the proof of Theorem 1 we obtain the following conjecture\\
\\
\textbf{Conjecture}\\
\emph{Let $u(x)$ be elliptic (meromorphic) function with finite number of poles in a fundamental parallelogram  ( in $\mathbb{C}$ and hasn't isolated pole at infinity ). Assume $S=\sum \limits_{i=1}^{i=m} n_i + 1$, where $m$ is number of poles. If $\varphi_{i,-4}=n_i(4n_i+1)(4n_i+3)(4n_i+4)$, $\varphi_{i,4k-l}=0$, $\varphi_{i, 4r-1}=\varphi_{i, 4r-3}=0$, where $k=0,...,n_i$, $r=n_i+1,...,S$ and $l=1,2,3$, then $L_4=\partial_x^4 + u(x)$ commutes with a differential operator of order $4S+2$}.
\\
\\
\textbf{Theorem 1.5}\\
\label{thm2} \emph{If $L_4 $ commutes with a differential operator $M$ of order $4g+2$, $M$ and $L_4$ are operators of rank 2 and $u(x)$ has pole at $a_i$, then solutions of the equation $\psi^{(4)}(x) + u(x)\psi(x)=\lambda \psi(x)$ have singularities at $a_i$ of the following type $x^{\sigma_{i,r}}g(x)$, where $r=1,2,3,4$ and $g(x)$ is holomorphic at $a_i$,\\
$\sigma_{i,1}=\frac{1}{2}(1-4n_i-\sqrt{1-16n_i-16n_i^2})$\\
$\sigma_{i,2}=\frac{1}{2}(1-4n_i+\sqrt{1-16n_i-16n_i^2})$\\
$\sigma_{i,3}=\frac{1}{2}(5+4n_i-\sqrt{1-16n_i-16n_i^2})$\\
$\sigma_{i,4}=\frac{1}{2}(5+4n_i+\sqrt{1-16n_i-16n_i^2})$.\\
This means that common eigenfunctions of commuting operators always have branch points. Hence $L_4$ doesn't commute with an operator of odd order.}\\
\\
\\
The author wishes to express gratitude to Professor O. I. Mokhov for advice and help in writing this paper.

\begin{center}
\textbf{2. Commuting differential operators of rank 2}
\end{center}

Consider the operator
\begin{equation}\label{operator}
L=(\partial^2_x + V(x))^2 + W(x)
\end{equation}
From \cite{Mironov} it follows that the operator $L$ commutes with an operator $M$ of order $4g+2$  and the spectral curve of $L$ and $M$ is hyperelliptic  curve of genus $g$ and hence operators $L$ and $M$ are operators of rank 2, if and only if there exists a polynomial
\begin{equation*}
Q=z^g + a_1(x)z^{g-1} + a_2(x)z^{g-2} + ...+a_{g-1}(x)z + a_g(x)
\end{equation*}
that the following relation is satisfied
\begin{equation}
Q^{(5)} + 4VQ''' + 6V'Q'' + 2Q'(2z-2W+V'') - 2QW'\equiv 0,
\end{equation}
$Q'$ means $\partial_xQ$. The spectral curve has the form
\begin{equation}
4w^2=4F(z)=4(z-W)Q^2 - 4V(Q')^2 + (Q'')^2 - 2Q'Q''' + 2Q(2V'Q' + 4VQ'' + Q^{(4)}).
\end{equation}
If $V(x)\equiv 0$, then we have
\begin{equation}\label{rec}
Q^{(5)}  + 2Q'(2z-2W) - 2QW'\equiv 0,
\end{equation}
\begin{equation*}
4w^2=4F(z)= 4zQ^2 + (Q'')^2 - 2Q'Q''' + 2Q Q^{(4)}.
\end{equation*}
Using (\ref{rec}) we get
\begin{equation*}
4Q'z\equiv -Q^{(5)}  + 2QW' + 4Q'W.
\end{equation*}
So we have the following system :
\begin{equation*}
\left\{
 \begin{array}{l}
a_1=W/2 + C_1\\
4a_2' = -a_1^{(5)}  + 2a_1W' + 4a_1'W\\
...\\
4a_{i+1}' = -a_i^{(5)}  + 2a_iW' + 4a_i'W\\
...\\
4a_{g}' = -a_{g-1}^{(5)}  + 2a_{g-1}W' + 4a_{g-1}'W\\
0 = -a_g^{(5)}  + 2a_gW' + 4a_g'W
 \end{array}
\right.  \label{sys}
\end{equation*}
It is clear from the system that the result in \cite{Mironov} can be formulated in the following way. Let $a_1=W/2 + C_1$, where $C_1$ is an arbitrary constant. Define $a_i$ by recursion
\begin{equation}
a_{i+1}=C_{i+1} + \frac{1}{4}\int(-a_i^{(5)} + 2a_iW' + 4a_i'W)dx .
\end{equation}
We see that if $V(x)\equiv 0$, then (\ref{operator}) commutes with an operator of order $4g+2$ and these operators are operators of rank 2 if and only if $a_{g+1}\equiv const$. \\
From the proof of Theorem 2 in \cite{Mironov} easy to see that $a_i(x)$ is a polynomial in $W(x)$, $W'(x)$, $W''(x)$,... But we will prove this fact by direct calculations.  We have
\begin{equation*}
\begin{gathered}
\int(W(x)a_1'(x))dx =\int(\frac{W(x)W'(x)}{2}) = \frac{W^2(x)}{4},\\
\int(3W^2(x)a_1'(x) - a_1^{(5)}(x)W(x))dx =\\
 =\frac{1}{2}W^3(x) - \frac{1}{4}(W''(x))^2 + \frac{1}{2}W'(x)W'''(x) - \frac{1}{2}W(x)W^{(4)}(x).
\end{gathered}
\end{equation*}
Let us suppose that $a_k$, $\int W(x)a_k'(x)dx$ and $\int(3W^2(x)a_k'(x) - a_k^{(5)}(x)W(x))dx$ are polynomials in $W(x)$, $W'(x)$,$W''(x)$,... We obtain
\begin{equation*}
\begin{gathered}
a_{k+1}(x) =  C_{k+1} + \frac{1}{4}\int(-a_k^{(5)}(x) + 2a_k(x)W'(x) + 4a_k'(x)W(x))dx =\\
=C_{k+1} - \frac{1}{4}\int a_k^{(5)}(x)dx + \frac{1}{2}\int \frac{d(a_k(x)W(x))}{dx}dx + \frac{1}{2}\int a_k'(x)W(x)dx.
\end{gathered}
\end{equation*}
So, we see that $a_{k+1}(x)$ is polynomial in $W(x)$, $W'(x)$,$W''(x)$... And analogously easy to check that $\int W(x)a_{k+1}'(x)dx$ and $\int(3W^2(x)a_{k+1}'(x) - a_{k+1}^{(5)}(x)W(x))dx$ are polynomial  in $W(x)$, $W'(x)$,$W''(x)$...
\\
\\
Let $L_4$ be the operator
\begin{equation*}
L_4=\partial_x^4 + u(x).
\end{equation*}
The commutativity condition  of $L_4$ and $L_6$ is equivalent to the equation
\begin{equation*}
4C_1u'(x) + 6u(x)u'(x) - u^{(5)}(x)=0,
\end{equation*}
for some constant $C_1$. Consider differential equations in $u(x)$
\begin{equation}
0=f_{j+1}(x)=C_{i+1} + \int (u(x) f'_j(x) + \frac{u'(x) f_j(x)}{2} - \frac{f_j^{(5)}(x)}{4}) dx,
\end{equation}
where
\begin{equation*}
f_1=C_1 + \frac{u(x)}{2}.
\end{equation*}
We see that $L_4$ commutes with an operator of order $4g+2$ if and only if there exist constants $C_1,...,C_g$ such that $f_{g+1}\equiv const$.
\\
\begin{center}
\textbf{3. Proof of Theorem 1.1 and corollaries}
\end{center}
\begin{flushleft}
\textbf{3.1 Proof of theorem 1.1}
\end{flushleft}
Assume that $u(x)$  has poles at points $a_1,a_2,...$ in $\mathbb{C}$.\\
\\
\textbf{Lemma 3.1}\\
\emph{If $L_4$ commutes with an operator of order  $4g+2$, then $u(x)$ can have poles in $\mathbb{C}$ only of order 4.}
\\
\begin{proof}
Assume that $u(x)$ has pole of order $k$  at point $a_i$. Since $L_4$ commutes with an operator of order $4g+2$, we have $f'_{g+1} \equiv 0$. Suppose $f_g$ has pole at point $a_i$ of order $m$; then $u(x)f'_g(x)$ + $\frac{u'(x)f_g(x)}{2}$ has pole of order $k+m+1$ and $f^{(5)}_g(x)$ has pole of order $m+5$. Hence, we see that if  $f_{g+1} \equiv 0$, then k+m+1=m+5. We get $k=4$.\\
\end{proof}

By $A^k_{i,m}$ denote a coefficient in the term $(x-a_i)^m$ in Laurent series  of $f_k$ at point $a_i$. We see from Lemma 3.1 that $f_k$ has pole of order $4k$ at point $a_i$. Assume that Laurent series of $u(x)$ has the form\\
\begin{equation*}
u(x) = \frac{\varphi_{i,-4}}{(x-a_i)^4} + \frac{\varphi_{i,-3}}{(x-a_i)^3} + \frac{\varphi_{i,-2}}{(x-a_i)^2} + \frac{\varphi_{i,-1}}{x-a_i} +
\end{equation*}
\begin{equation*}
+\varphi_{i,0} + \varphi_{i,1}(x-a_i) + \varphi_{i,2}(x-a_i)^2 + \varphi_{i,3}(x-a_i)^3 + \varphi_{i,4}(x-a_i)^4+ + O((x-a_i)^5).
\end{equation*}
\textbf{Lemma 3.2}\\
\emph{If $L_4$ commutes with an operator of order $4g+2,$ then\\
$\varphi_{i,-4}=n_i(4n_i +1)(4n_i+3)(4n_i + 4)$, where $n_i \in \mathbb{N}$ and}
\begin{equation*}
A^{k+1}_{i, -4k-4}=\frac{(2k+1)A^k_{i,-4k}(\varphi_{-4} - k(4k+1)(4k+3)(4k+4))}{2k+2}.
\end{equation*}

\begin{proof}
In a neighborhood of $a_i$
\begin{equation*}
f_j=\frac{A_{i, -4j}^j}{(x-a_i)^{4j}} + O((x-a_i)^{-4j+1}).
\end{equation*}
Easy to see that
\begin{equation*}
f'_{j+1}=(\frac{\varphi_{i,-4}}{(x-a_i)^4} + O(\frac{1}{(x-a_i)^3}))(-\frac{4jA_{i, -4j}^j}{(x-a_i)^{4j+1}} + O((x-a_i)^{-4j}) +
\end{equation*}
\begin{equation*}
+(-\frac{2\varphi_{i,-4}}{(x-a_i)^5} + O(\frac{1}{(x-a_i)^4}))(\frac{A_{i, -4j}^j}{(x-a_i)^{4j}} + O((x-a_i)^{-4j+1}) +
\end{equation*}
\begin{equation*}
+\frac{j(4j+1)(4j+2)(4j+3)(4j+4)A_{i,-4j}^j}{(x-a_i)^{4j+5}} + O((x-a_i)^{-4j-4}) =
\end{equation*}
\begin{equation*}
= -\frac{(4j+2)(\varphi_{i,-4} - j(4j+1)(4j+3)(4j+4))A_{i, -4j}^j}{(x-a_i)^{4j+5}} + O((x-a_i)^{-4j-4}).
\end{equation*}
So, integrating we obtain that
\begin{equation*}
A^{j+1}_{i, -4j-4} = \frac{(2j+1)(\varphi_{i,-4} - j(4j+1)(4j+3)(4j+4))A_{i, -4j}^j}{2j+2}
\end{equation*}
If  $f'_{g+1}\equiv 0$, then the leading term of principal part must vanish. Hence, $\varphi_{i,-4}=n_i(4n_i +1)(4n_i+3)(4n_i+4)$ for some $n_i$.\\
\end{proof}
\begin{flushleft}
\textbf{Lemma 3.3}\\
\end{flushleft}
\emph{Suppose $j < n_i+1$ and in a neighborhood of $a_i$}
\begin{equation*}
f_j=\frac{A_{i, -4j}^j}{(x-a_i)^{4j}} + O((x-a_i)^{-4j+1});
\end{equation*}
\emph{then $A_{i, -4j}^j > 0$ and if $j \geqslant n_i+1$, then  $A^{j}_{i, -4j} = 0.$}
\begin{proof}
Since $\varphi_{i,-4}=n_i(4n_i+1)(4n_i+3)(4n_i+4)$ and
\begin{equation}
A_{i, -4j-4}^{j+1} = \frac{(2j+1)(\varphi_{i,-4} - j(4j+1)(4j+3)(4j+4))A_{i,4j}^{j}}{2j+2},
\end{equation}
we see that if $j < n_i+1$, then $A_{i, -4j}^j > 0$.\\
If $j \geqslant n_i+1$, then $A_{i, -4j}^{j}=0$.
\end{proof}
Let us prove that $u(x)$ can't have isolated pole at infinity. Assume the converse. Then $u(x)$ has pole at infinity of order $m$. We have
\begin{equation*}
f_1 = A^1_{\infty, m}x^{m} + O(x^{m-1}) =\frac{\varphi_{\infty, m}}{2} + O(x^{m-1})
\end{equation*}
But $f_2$ has pole at infinity of order $2m$ and $A^2_{\infty, 2m} = \frac{2mA^1_{\infty, m}\varphi_{\infty, m} + m\varphi_{\infty, m}A^1_{\infty, m}}{4m} \neq 0.$ In general $A^{k+1}_{\infty, (k+1)m} = \frac{2kmA^k_{\infty, km}\varphi_{\infty, m} + m\varphi_{\infty, m}A^k_{\infty, km}}{2m(k+1)}=\frac{\varphi_{\infty, m}A^k_{\infty, mk}(2k+1)}{2(k+1)} \neq 0$. So, order of pole increases and there is no $k$ such that $f_k$ vanishes.\\
\\
\textbf{Lemma 3.4}\\
\emph{If $L_4$ commutes with an operator of order $4g+2$, then\\ $\varphi_{i,-3}=\varphi_{i,-2}=\varphi_{i,-1}=0$ for all $i$.}\\
\begin{proof}
By definition, $f_1(x) = C_1 + \frac{u(x)}{2}$ in a neighborhood of $a_i$ has the form\\
\begin{equation*}
f_1=\frac{\varphi_{i,-4}}{2(x-a_i)^4} + \frac{\varphi_{i,-l}}{2(x-a_i)^l} + O((x-a_i)^{-l+1}),
\end{equation*}
where $l=1,2,3$. Therefore, $f'_2(x)=u(x)f'_1(x) + \frac{u'(x)f_1(x)}{2} - \frac{f_1^{(5)}(x)}{4}$ in a neighborhood of $a_i$
\begin{equation*}
(\frac{\varphi_{i,-4}}{(x-a_i)^4} + \frac{\varphi_{i,-l}}{(x-a_i)^l} + O((x-a_i)^{-l+1}))(-\frac{2\varphi_{i,-4}}{(x-a_i)^5} - \frac{l\varphi_{i,-l}}{2(x-a_i)^{l+1}} + O((x-a_i)^{-l})) +
\end{equation*}
\begin{equation*}
+(-\frac{\varphi_{i,-4}}{(x-a_i)^5} -\frac{l\varphi_{i,-l}}{4(x-a_i)^{l+1}} + O((x-a_i)^{-l}))(\frac{\varphi_{i,-4}}{(x-a_i)^4} + \frac{\varphi_{i,-l}}{(x-a_i)^l} + O((x-a_i)^{-l+1})) +
\end{equation*}
\begin{equation*}
+\frac{840\varphi_{i,-4}}{(x-a_i)^{9}} + \frac{l(l+1)(l+2)(l+3)(l+4)\varphi_{i,-l}}{8(x-a_i)^{l+5}} + O((x-a_i)^{-l-4}) =
\end{equation*}
\begin{equation*}
= -\frac{3\varphi_{i,-4}(\varphi_{i,-4}-280)}{(x-a_i)^{9}} - \frac{\varphi_{i,-l}(4+l)(6\varphi_{i,-4} -l(l+1)(l+2)(l+3))}{8(x-a_i)^{l+5}}  + O((x-a_i)^{-l-4})
\end{equation*}
Integrating, we obtain
\begin{equation*}
f_2(x) = \frac{3\varphi_{i,-4}(\varphi_{i,-4}-280)}{8(x-a_i)^{8}} + \frac{\varphi_{i,-l}(6\varphi_{i,-4} -l(l+1)(l+2)(l+3))}{8(x-a_i)^{l+4}}  + O((x-a_i)^{-l-3})
\end{equation*}
Note that $6\varphi_{i,-4} -l(l+1)(l+2)(l+3) > 0$.  Consider $f_k$, where $k\leqslant n_i$.
\begin{equation*}
f_k=\frac{A_{i, -4k}^k}{(x-a_i)^{4k}} + \frac{A_{i, -4k+4-l}^k}{(x-a_i)^{4k-4+l}} + O((x-a_i)^{-4k+3-l}).
\end{equation*}
Let us prove that
\begin{equation}
A_{i, -4k-l}^{k+1}=\varphi_{i,-l}K_{i, -4k-l}^{k+1} \quad for \quad k \leqslant n_i,
\end{equation}
where $K_{i, -4k-l}^{k+1}>0$ and doesn't depend on $\varphi_{i,-l}.$ The proof is by induction on $k$. We checked this for $k=2$. By the induction hypothesis, $A_{i,-4k+4 -l}^k=\varphi_{i,-l}K_{i, -4k+4-l}^k$ and $K_{i, -4k+4+l}^k >0$.\\
Easy to see that\\
\begin{equation*}
f_{k+1}= \frac{(2k+1)(\varphi_{i,-4} - k(4k+1)(4k+3)(4k+4))A_{i, -4k}^{k}}{(2k+2)(x-a_i)^{4k+4}} +
\end{equation*}
\begin{equation*}
+\frac{\varphi_{i,-l}A^k_{i,-4k}(8k+l)}{2(4k+l)(x-a_i)^{4k+l}} +
\end{equation*}
\begin{equation*}
\frac{(4k-2+l)(4\varphi_{i,-4} - (4k-4+l)(4k-3+l)(4k-1+l)(4k+l))A^k_{i,-4k+4-l} }{4(4k+l)(x-a_i)^{4k+l}} + ...
\end{equation*}
\begin{equation*}
=\frac{A_{i, -4k-4}^{k+1}}{(x-a_i)^{4k+4}} + \frac{A_{i, -4k-l}^{k+1}}{(x-a_i)^{4k+l}} + O((x-a_i)^{-4k+1-l})
\end{equation*}
\\
We see that $A^{k+1}_{i,-4k-l}=\varphi_{i,-l}K^{k+1}_{i,-4k-l}$. We obtain that $K^{k+1}_{i,-4k-l}>0$ because $A_{i,-4k}^{k}>0$, where $k \leqslant n_i$. So, (8) is proved.     \\
We know from Lemma 3.2 that $A_{i, -4n_i - 4}^{n_i + 1}=0$. If $k\geqslant n_i + 1$, then
\begin{equation*}
A^{k+1}_{i,-4k-l} =\frac{A_{i,-4k + 4 -l}^k(4k-2+l)(4\varphi_{i,-4} -  (4k-4+l)(4k-3+l)(4k-1+l)(4k+l))}{4(4k+l)}.
\end{equation*}
Easy to see that if $k\geqslant n_i + 1$, then
\begin{equation*}
4n_i(4n_i+1)(4n_i+3)(4n_i + 4) - (4k-4+l)(4k-3+l)(4k-1+l)(4k+l) ) < 0.
\end{equation*}
Finally we obtain that if there exists a  $g$ such that $f_{g+1} \equiv 0$, then $A^{k}_{i,-4k+4-l}=\varphi_{i,-l}K^{k}_{i,-4k+4-l}=0$ for some $k$. But $K^{k}_{i,-4k+4-l} \neq 0$ for all $k$. Hence, $\varphi_{i,-l}=0$.
\end{proof}
\begin{flushleft}
Now let us prove the main part of Theorem 1.1.
\end{flushleft}
From Lemma 3.4 we know that  $f_1(x) = C_1 + \frac{u(x)}{2}$  in a neighborhood of $a_i$ has the form\\
\begin{equation*}
f_1=\frac{\varphi_{i,-4}}{2(x-a_i)^4} + \frac{\varphi_{i,0}}{2} + \widetilde{C}_1 + \frac{\varphi_{i,4-l}(x-a_i)^{4-l}}{2} + O((x-a_i)^{5-l}),
\end{equation*}
where $l=1,2,3$. Then $f'_2(x)=u(x)f'_1(x) + \frac{u'(x)f_1(x)}{2} - \frac{f_1^{(5)}(x)}{4}$  in a neighborhood of $a_i$
\begin{equation}
f_2= \frac{3(\varphi_{i,-4} - 280)\varphi_{i,-4}}{8(x-a_i)^8} + \frac{\varphi_{i,-4}(3\varphi_{i,0} + 2C_1)}{4(x-a_i)^4} + \frac{3\varphi_{i,-4}\varphi_{i,4-l}}{4(x-a_i)^{l}} + O((x-a_i)^{-l+1}).
\end{equation}
We see that the coefficient in the term  $(x-a_i)^{-l}$ equals $\varphi_{i,l}\frac{3\varphi_{i,-4}}{4}$. Consider $f_k$
\begin{equation*}
f_k= \frac{A^k_{i, -4k}}{(x-a_i)^{4k}} + \frac{A^k_{i,-4k+4}}{(x-a_i)^{4k-4}} + \frac{A^k_{i,-4k+8-l}}{(x-a_i)^{4k-8+l}} + O((x-a_i)^{-4k+9-l})
\end{equation*}
Let us show that $A^k_{i, -4k+8-l}=\varphi_{i,l}K^k_{i, -4k+8-l}$, where $K^k_{i, -4k+8-l} \neq 0$ and doesn't depend on $\varphi_{i,4-l}$. The proof is by induction on $k$. We checked this for $k=2$.
\begin{equation}
\begin{gathered}
f_{k+1}= \frac{A^{k+1}_{i,-4k-4}}{(x-a_i)^{4k+4}} + \frac{A^{k+1}_{i,-4k}}{(x-a_i)^{4k}} + \frac{(8k-4+l)\varphi_{i, 4-l}A^k_{i, -4k}}{2(4k-4+l)(x-a_i)^{4k-4+l}} +\\
+\frac{(4k-6+l)(4\varphi_{i, -4} - (4k-8+l)(4k-7+l)(4k-5+l)(4k-4+l))A^k_{i, -4k+8-l}}{4(4k-4+l)(x-a_i)^{4k-4+l}}+\\
 + O((x-a_i)^{-4k+5-l})
\end{gathered}
\end{equation}
So, we have
\begin{equation*}
\begin{gathered}
K^{k+1}_{i, -4k+4-l}=\frac{(8k-4+l)A^k_{i, -4k}}{2(4k-4+l)} + \\
+\frac{(4k-6+l)(4\varphi_{i, -4} - (4k-8+l)(4k-7+l)(4k-5+l)(4k-4+l))K^k_{i, -4k+8-l}}{4(4k-4+l)}.
\end{gathered}
\end{equation*}
We see that  if $k \leqslant n_i+1$, then $K^{k+1}_{i, -4k+4-l} > 0$. If $k>n_i+1$, then  from Lemma 3.3 we obtain that  $A^k_{i, -4k}=0$. Hence, if $k>n_i+1$, then $K^k_{i, -4k+8-l} <0.$ But $A^k_{i, -4k+8-l}=\varphi_{i,l}K^k_{i, -4k+8-l}$ and we must find a $g$ such that $f_{g+1} \equiv 0$. It now follows that  $A^k_{i, -4k+8-l}=0$ for some $k$ $\Leftrightarrow \varphi_{i,l}$ = 0.\\
\\
In general, assume $\varphi_{i, 4k-l}=0$, where $k=1,...,m-1$, $l=1,2,3$ and $m \leqslant n_i$. We have already checked this for $m=2$. Let us prove that $\varphi_{i, 4m-l}=0$. We have
\begin{equation}
f_1=\frac{\varphi_{i, -4}}{2(x-a_i)^4} + \frac{\varphi_{i, 0}}{2} + C_1 + \sum_{t=1}^{t=m-1}\frac{\varphi_{i, 4t}(x-a_i)^{4t}}{2} + \frac{\varphi_{i, 4m-l}(x-a_i)^{4m-l}}{2} + ...
\end{equation}
Then
\begin{equation*}
f_k = \sum^{t=m-k}_{t=-k} A^k_{i, 4t}(x-a_i)^{4t} + A^k_{i,4m-4k+4-l}(x-a_i)^{4m-4k+4-l} +...
\end{equation*}
\begin{equation*}
f_{k+1}  = \sum^{t=m-k-1}_{t=-k-1} A^{k+1}_{i, 4t}(x-a_i)^{4t}  + A^{k+1}_{i,4m-4k-l}(x-a_i)^{4m-4k-l} + ...
\end{equation*}
Calculations show that
\begin{equation}
\begin{gathered}
 A^{k+1}_{i, 4m-4k-l} = \\
= \frac{(4m-4k+2-l) A^{k}_{i, 4m-4k+4-l}}{4(4m-4k-l)} \times\\
\times (4\varphi_{i, -4} - (4m-4k+4-l)(4m-4k+3-l)(4m-4k+1-l)(4m-4k-l)) +
\end{gathered}
\end{equation}
\begin{equation*}
 +  \frac{(4m-8k-l)\varphi_{i, 4m-l}A^k_{i,-4k}}{2(4m-4k-l)}.
\end{equation*}
If we take $m=1$ in (12), we get (10).\\
Since $A^1_{i, 4m-l}=\frac{\varphi_{i, 4m-l}}{2}$, we obtain from (12) that $A^{k+1}_{i,4m-4k-l}=\varphi_{i,4m-l}K^{k+1}_{i,4m-4k-l}$, where $K^{k+1}_{i,4m-4k-l}$ doesn't depend on $\varphi_{i, 4m-l}.$
\begin{equation}
\begin{gathered}
 K^{k+1}_{i, 4m-4k-l} = \\
= \frac{(4m-4k+2-l) K^{k}_{i, 4m-4k+4-l}}{4(4m-4k-l)} \times \\
\times (4\varphi_{i, -4} - (4m-4k+4-l)(4m-4k+3-l)(4m-4k+1-l)(4m-4k-l)) +
\end{gathered}
\end{equation}
\begin{equation*}
 +  \frac{(4m-8k-l)A^k_{i,-4k}}{2(4m-4k-l)}.
\end{equation*}

\begin{flushleft}
\textbf{Lemma 3.5}
\end{flushleft}
\emph{Number $K^{n_i+1}_{i,4m -4k-l} \neq 0$ for  $m\leqslant n_i$.}
\begin{proof}
Consider two cases.
\\
$1)$ $l=2$.\\
From (13) we obtain that
\begin{equation*}
 A^{m+1}_{i, -2} =  \frac{(4m+2)\varphi_{i, 4m-2}A^m_{i,-4m}}{4}=\varphi_{i, 4m-2}\frac{(4m+2)A^m_{i,-4m}}{4}=\varphi_{i, 4m-l}K^{m+1}_{i,-2}.
\end{equation*}
We see from Lemma 3.2 that if $m \leqslant n_i$, then $K^{m+1}_{i,-2} >0$. But for $k > m$ expression (13) is positive. So $K^{k+1}_{i, 4m-4k-2}>0$ for all $k>m$. \\
\\
$2)$ $l=1$ or $l=3$.
\\
If $K^{n_i+1}_{i, 4m-4n_i-l}=0$, then from (13) we get
\begin{equation*}
\begin{gathered}
\frac{(4m-8n_i-l)A^{n_i}_{i,-4n_i}}{2(4m-4n_i-l)}=\\
=-\frac{(4m-4n_i+2-l) A^{n_i}_{i, 4m-4n_i+4-l}}{4(4m-4n_i-l)} \times \\
\times (4\varphi_{i, -4} - (4m-4n_i+4-l)(4m-4n_i+3-l)(4m-4n_i+1-l)(4m-4n_i-l)).
\end{gathered}
\end{equation*}
But calculating $A^{n_i}_{i, 4m-4n_i+4-l}$ recursively using (12) we see that the expression above is not true.\\
\end{proof}
We know from Lemma 3.3 that $A^k_{i,-4k} = 0$ for all $k \geqslant n_i+1$. So for $k \geqslant n_i +1$ we have
\begin{equation*}
\begin{gathered}
K^{k+1}_{i, 4m-4k-l} = \\
=\frac{(4m-4k+2-l) K^{k}_{i, 4m-4k+4-l}}{4(4m-4k-l)} \times \\
\times (4\varphi_{i, -4} - (4m-4k+4-l)(4m-4k+3-l)(4m-4k+1-l)(4m-4k-l)) \neq 0
\end{gathered}
\end{equation*}
So, if there exists a $g$ such that $f_{g+1}=0$, then for some $k$ coefficient $A^{k+1}_{i, 4m-4k-l}=\varphi_{i, 4m-l}K^k_{i, 4m-4k-l}=0 \Leftrightarrow \varphi_{i, 4m-l}=0$.\\
\\
\textbf{Now suppose that $g>n_i$}. Then
\begin{equation*}
f_1=\frac{\varphi_{i, -4}}{2(x-a_i)^4} + \frac{\varphi_{i, 0}}{2} + C_1 + \sum_{t=1}^{t=m-1}\frac{\varphi_{i, 4t}(x-a_i)^{4t}}{2} + \frac{\varphi_{i, 4m-l}(x-a_i)^{4m-l}}{2} + ...,
\end{equation*}
where $m \geqslant n_i+1$. We see that
\begin{equation*}
f_k = \sum^{t=m-k}_{t=-k} A^k_{i, 4t}(x-a_i)^{4t} + A^k_{i,4m-4k+4-l}(x-a_i)^{4m-4k+4-l} +...
\end{equation*}
If $k>n_i$, we get
\begin{equation*}
\begin{gathered}
K^{k+1}_{i, 4m-4k-l} =\\
\frac{(4m-4k+2-l) K^{k}_{i, 4m-4k+4-l}}{4(4m-4k-l)} \times \\
\times (4\varphi_{i, -4} - (4m-4k+4-l)(4m-4k+3-l)(4m-4k+1-l)(4m-4k-l)).
\end{gathered}
\end{equation*}
So, we obtain $K^{k+1}_{i,4m-4k-l}  \neq 0$ for $l=1,3$ and $k>n_i$  because
\begin{equation*}
4\varphi_{i, -4} - (4m-4k+4-l)(4m-4k+3-l)(4m-4k+1-l)(4m-4k-l)) \neq 0,
\end{equation*}
If $k=m>n_i$ and $l=2$, then $K^{m+1}_{i, -2}=0$ because of factor $(4m-4k+2-l)$.  But if $L_4$ commutes with an operator of order  $4g+2$, then there exists a $k\leqslant g+1$ such that $A^k_{i, 4m-4k+3}=A^k_{i, 4m-4k+1}=0 \Leftrightarrow \varphi_{i, 4m-1}=\varphi_{i, 4m-3}=0$.\\
\\
\textbf{Theorem 1.1 is proved.}
\begin{flushleft}
\textbf{3.2 Proof of corollaries and examples}
\end{flushleft}
We know that
\begin{equation*}
f_2= \frac{3(\varphi_{i,-4} - 280)\varphi_{i,-4}}{8(x-a_i)^8} + \frac{\varphi_{i,-4}(3\varphi_{i,0} + \widetilde{C}_1)}{4(x-a_i)^4} + \widetilde{C}_2 + O((x-a_i)),
\end{equation*}
where $\widetilde{C}_i$ are constants and depend on $C_i$. Also $A^{n_i+1}_{i,-4(n_i+1)}=0$ and that's why
\begin{equation*}
f_{n_i+1}= \sum\limits_{t=n_i}^{t=1}\frac{A^{n_i+1}_{i, -4t}}{(x-a_i)^{4t}} + \widetilde{C}_{n_i+1} + O((x-a_i)),
\end{equation*}
We must choose constants $C_1,...C_{g}$ to vanish principal part of Laurent series of function $f_{n_i+1}$.  This is always possible to do because $A^{n_i+1}_{i,-4n_i}$ linearly depend on $C_1$, $A^{n_i+1}_{i,-4n_i+4}$ linearly depend on $C_1$ and $C_2$, $A^{n_i+1}_{i,-4n_i+8}$ depend on $C_1, C_2, C_3$,  $A^{n_i+1}_{i,-4}$ depend linearly on $C_1,...,C_{n_i}$. We know that only constant function can be holomorphic and bounded function on $\mathbb{C}$. Hence  $f'_{n_i+1} \equiv 0$.\\
Corollary 1.2 is proved.\\
\\
If we take $C_1=-42g_2$, then we obtain Example 1.4.

\begin{flushleft}
\textbf{3.3 Conjecture}
\end{flushleft}
Consider $f_k$, where $k \geqslant n_i + 1$. Assume that $u(x)$ has $m$ poles (in the fundamental parallelogram or $\mathbb{C}$). Take $S=\sum\limits^{i=m}_{i=1} n_i +1 $. If $\varphi_{i,4k-l}=0$, $\varphi_{i, 4r-1}=\varphi_{i, 4r-3}=0$, where $k=0,...,n_i$ and $r=n_i+1,...,S$. We see from Lemma 3.3 that degree of pole at $a_i$ is not greater than $n_i$ for any $k$ . We get that $A^k_{i, -4n_i}$ linearly depends on $C_1,...,C_{k-n_i}$, $A^k_{i, -4n_i+4}$ linearly depends on $C_1,...,C_{k-n_i+1}$, $A^k_{i,-4}$ depends on $C_1,...,C_k$.   So, we must choose constants to vanish $S-1$ terms and we have $S$ constants. Apparently this linear system always has solutions , but it is not proved.

\begin{center}
\textbf{4. Proof of Theorem 1.5}
\end{center}
Consider the differential equation
\begin{equation}
\frac{d^nw(x)}{dx^n} + \widetilde{P}_{1}(x)\frac{d^{n-1}w(x)}{dx^{n-1}}+...+\widetilde{P}_{n-1}(x)\frac{dw(x)}{dx}+\widetilde{P}_n(x)w(x)=0.
\end{equation}
We know from theory of ordinary differential equations (see ~\cite{Ince} chapter 16) that if $w(x)$ has singularity at point $a$, then $P_i(x)$ has singularity at point $a$ for some $i$. A singularity is called regular if coefficients $P_k$ have pole of order not greater than $k$ for all $k$. Without loss of generality it can be assumed that $a=0$ and we can write (14) in the following form
\begin{equation*}
x^n\frac{d^nw(x)}{dx^n} + x^{n-1}P_{1}(x)\frac{d^{n-1}w(x)}{dx^{n-1}}+...+xP_{n-1}(x)\frac{dw(x)}{dx}+P_n(x)w(x)=0,
\end{equation*}
where coefficients $P_i(x)$ haven't pole at point $a$.
. Solutions in a neighborhood of regular singularity  have the form (see ~\cite{Ince} chapter 16)
\begin{equation}
w(x)=\sum\limits^{\infty}_{m=0}c_m x^{m+\sigma}.
\end{equation}
By  $L$ denote the operator $x^n\frac{d^nw(x)}{dx^n} + x^{n-1}P_{1}(x)\frac{d^{n-1}w(x)}{dx^{n-1}}+...+xP_{n-1}(x)\frac{dw(x)}{dx}+P_n(x)w(x)=0$ and\\ $[\sigma+m]_n=(\sigma+m)(\sigma+m-1)...(\sigma+n-m+1)$. We get
\begin{equation*}
Lw=L(\sum\limits^{\infty}_{m=0}c_m x^{m+\sigma}) = \sum\limits^{\infty}_{m=0} c_mx^{m+\sigma}f(x,m+\sigma),
\end{equation*}
where $f(x, m+\sigma)=[\sigma+m]_n + P_1(x)[\sigma+m]_{n-1} + [\sigma+m]_1P_{n-1}(x) + P_n(x)=\sum\limits^{\infty}_{\lambda=0}f_{\lambda}(m+\sigma)x^{\lambda}$. If $Lw=0$, then\\
\\
$c_0f_0(\sigma)=0$\\
$c_1f_0(\sigma +1) + c_0f_1(\sigma)=0$\\
................................................\\
$c_mf_0(\sigma + m) + c_{m-1}f_1(\sigma+m-1)+...+c_0f_m(\sigma)=0$.\\
\\
Since $c_0\neq 0$, it follows that\\
\begin{equation}
f_0(\sigma)=[\sigma]_n + [\sigma]_{n-1}P_1(0)+...+[\sigma]_1P_{n-1}(0) + P_n(0)=0.
\end{equation}
Let $\sigma$ be  a root of (16). If $f_0(\sigma + m)\neq 0$, where $m$ is integer, then we can find constants $c_m$ by the following formulas (see ~\cite{Ince} chapter 16)
\begin{equation*}
c_m=\frac{(-1)^mc_0F_m(\sigma)}{f_0(\sigma+1)f_0(\sigma+2)...f_0(\sigma+m)},
\end{equation*}
where
$$F_m(\sigma)=\begin{vmatrix}
f_1(\sigma+m-1), & f_2(\sigma+m-2), & ..., & f_{m-1}(\sigma+1), & f_m(\sigma)  \\
f_0(\sigma+m-1), & f_1(\sigma+m-2), & ..., & f_{m-2}(\sigma+1), & f_{m-1}(\sigma) \\
0,               & f_0(\sigma+m-2), & ..., & f_{m-3}(\sigma+1), & f_{m-2}(\sigma) \\
......           & ......           & ..., & .......,           & ......\\
0,               & 0,               & ..., & f_0(\sigma+1),     & f_1(\sigma)
\end{vmatrix}$$
So, if $\sigma_i$, $\sigma_j$ are roots of $f_0(\sigma)=0$ and $\sigma_i - \sigma_j$ is not integer for all $i\neq j$, then all solutions of (14) have the form $x^{\sigma}g(x)$, where  $g(x)$ is a holomorphic function. If for some $i,j$ number $\sigma_i - \sigma_j \in \mathbb{Z}$, then solutions have the form $g_0(x)ln^kx + g_1(x)ln^{k-1}x + ...+g_k(x)$, where $g_i(x)$ could be equal zero.
\\
Consider the operator $L_4=\partial^4_x + u(x)$. We proved before that if $L_4$ commutes with an operator of order $4g+2$, then $u(x)$ can have pole only of order $4$ and $\varphi_{i,-4}=n_i(4n_i+1)(4n_i+3)(4n_i+4)$. This means that eigenfunctions of $L_4$ have regular singularities. Solutions of (16) have the form\\
$\sigma_{i,1}=\frac{1}{2}(1-4n_i-\sqrt{1-16n_i-16n_i^2})$\\
$\sigma_{i,2}=\frac{1}{2}(1-4n_i+\sqrt{1-16n_i-16n_i^2})$\\
$\sigma_{i,3}=\frac{1}{2}(5+4n_i-\sqrt{1-16n_i-16n_i^2})$\\
$\sigma_{i,4}=\frac{1}{2}(5+4n_i+\sqrt{1-16n_i-16n_i^2})$.\\
This means that solutions can have logarithmic terms because $\sigma_{i,3} - \sigma_{i,1}=\sigma_{i,4} - \sigma_{i,2} = 4n_i + 2$. But we proved  that if $L_4$ commutes with operator $4g+2$, then $\varphi_{4k-l}=0$, $k=0,...,n_i$, $l=1,2,3$. Hence, eigenfunctions in a neighborhood of pole $a_i$ have the form
\begin{equation*}
\psi_{i,k}(x)=\alpha_0x^{\sigma_{i,k}} + \alpha_4x^{\sigma_{i,k}+4} + ... + \alpha_{4n_i}x^{\sigma_{i,k} + 4n_i} + \alpha_{4n_i+1}x^{\sigma_{i,k} + 4n_i +1} + \alpha_{4n_i+2}x^{\sigma_{i,k} + 4n_i+2} + ...
\end{equation*}
Easy to check that we can find all coefficients explicitly. So logarithmic terms will not appear.

\end{document}